\documentclass[final]{siamltex}
\usepackage{amsmath,amssymb,graphics}

\newcommand{\R}{\mathbb{R}}
\renewcommand{\vec}[1]{\mathbf{#1}}
\newcommand{\eps}{\varepsilon}
\newcommand{\var}[1]{\ensuremath{\text{Var}(#1)}}
\newcommand{\T}{\mathbb{T}}
\newcommand{\N}{\mathbb{N}}
\newcommand{\erf}{\mathrm{erf}}  
\newcommand{\zero}{\{ \vec{0} \}}

\begin{document}

\begin{flushright}
{\it SIAM J. Appl. Math. \textbf{63}, 1615 (2003)}
\end{flushright}

\title{System of phase oscillators with diagonalizable interaction}
\author{Takashi Nishikawa\footnotemark[2]
\and Frank C. Hoppensteadt\footnotemark[2]}

\renewcommand{\thefootnote}{\fnsymbol{footnote}}
\footnotetext[2]{Department of Mathematics, Arizona State
University, Tempe, Arizona 85287-1804, USA (tnishi@chaos6.la.asu.edu, fchoppen@asu.edu).  The research of the first author was supported by DARPA/ONR grant N00014-01-1-0943.  The research of the second author was partially supported by NSF grant DMS-0109001.}
\renewcommand{\thefootnote}{\arabic{footnote}}

\maketitle

\begin{abstract}
We consider a system of $N$ phase oscillators having randomly
distributed natural frequencies and diagonalizable interactions among
the oscillators.  We show that, in the limit of $N \to \infty$, all
solutions of such a system are incoherent with probability one for any
strength of coupling, which implies that there is no sharp transition
from incoherence to coherence as the coupling strength is increased,
in striking contrast to Kuramoto's (special) oscillator system.
\end{abstract}

\begin{keywords}
Network of phase oscillators, Kuramoto model
\end{keywords}

\begin{AMS}
34C15, 
37N25, 
37N20 
\end{AMS}

\pagestyle{myheadings}
\thispagestyle{plain}
\markboth{TAKASHI NISHIKAWA AND FRANK C. HOPPENSTEADT}
{DIAGONALIZABLE SYSTEM OF PHASE OSCILLATORS}

\section{Introduction}
\label{sec:intro}

Synchronization of coupled oscillators is a ubiquitous phenomenon in
natural and artificial systems.  Examples include synchronization of
pacemaker cells of the heart~\cite{mechaels87,peskin75}, rhythmic
activities in the brain~\cite{cray94,singer95}, synchronous flashing
of fireflies~\cite{buck76,buck88}, arrays of
lasers~\cite{jiang93,kourtchatov95}, and superconducting Josephson
junctions~\cite{wiesenfeld96,wiesenfeld98}.  Characterization of the
phenomenon using mathematical models has been a topic of great
interest for researchers in various scientific and engineering
disciplines.

Wiener~\cite{wiener58,wiener61}, who recognized the ubiquity of
synchronization phenomena in the real world, made a first attempt at
characterization using the Fourier integrals.  A more successful
approach was taken by Winfree~\cite{winfree67}, who used a population
of interacting limit-cycle oscillators to describe synchronization
properties.  He realized that if the interactions among the
oscillators are weak and the oscillators are nearly identical, the
separation of fast and slow timescales leads to a reduced model that
can be expressed in terms solely of the phase of each oscillator.
Kuramoto~\cite{kuramoto84} put this idea on a firmer foundation by
employing a perturbation method to show that the reduced equation has
a universal form.  His analysis of this model in the case of
mean-field coupling kicked off an avalanche of theoretical
investigations of his model and its generalizations.

More generally and rigorously, if each oscillator has an exponentially
stable limit-cycle and interactions among them are weak, the reduced
phase equation can be shown (see Theorem
9.1~\cite[p.~253]{hoppensteadt97a}) to have the form
\begin{equation}
\dot{\theta} = \omega + \eps f(\theta),\quad \theta \in \T^N,
\label{eqn:phase}
\end{equation}
where $\omega \in \R^N$ is the vector of natural frequencies of the
oscillators that are coupled to one another through the interaction
function $f : \T^N \to \R^N$, and $\eps > 0$ represents the overall
strength of the coupling.  The universal form of the interaction
function, derived by Kuramoto~\cite{kuramoto84} under the additional
assumption that the oscillators are almost identical, corresponds to
the choice
\begin{equation}
\label{eqn:kuramoto}
f_i(\theta) = \sum_{j=1}^N h_{ij}(\theta_j - \theta_i),
\quad i = 1, \ldots, N,
\end{equation}
where $f(\theta) = (f_1(\theta), \ldots, f_N(\theta))^T$.  The
mean-field model that he studied results when $h_{ij}(x) = \sin(x)/N$ for
all $i,j$.

Let $\theta(t)$ be a solution of~\eqref{eqn:phase}.  The
oscillators $i$ and $j$ are said to be \emph{locked} if $\lim_{t \to
\infty} \theta_i(t)/\theta_j(t) = 1$.  The solution is said to be
\emph{coherent} if all pairs of oscillators are locked.  If none of
the oscillator pairs are locked, the solution is \emph{incoherent}.  A
solution that is neither coherent nor incoherent is called
\emph{partially coherent}.  The main conclusion of Kuramoto's
work~\cite{kuramoto84} on his mean-field model is that in the limit of
$N \to \infty$ there exists a critical coupling strength $\eps_c$
such that for $\eps < \eps_c$ the solution is incoherent, but for
$\eps > \eps_c$ partially coherent solutions appear, for which the
fraction of locked oscillator pairs is nonzero.  Although his result
was important, since this behavior closely resembles the phase
transition phenomena widely observed in statistical physics, his
analysis is heuristic and makes assumptions about the symmetry of the
distribution of natural frequencies, which might not be necessary for
the results~\cite{strogatz93}.

In this paper, we consider a class of \emph{diagonalizable}
interaction functions, in which separation of variables is possible
after an appropriate coordinate transformation.  This allows us to
prove rigorously that for the system~\eqref{eqn:phase} with a generic
diagonalizable interaction function, if the solution is partially
coherent, then it is almost surely coherent.  This, together with the
fact that the probability of having a coherent solution goes to zero
in the limit of $N \to \infty$, leads to our main conclusion.  Namely,
for any $\eps > 0$, the solution is almost surely incoherent in the
limit of $N \to \infty$.  Our result shows that a diagonalizable
system of phase oscillators cannot exhibit a sudden transition from
incoherence to coherence, in sharp contrast to the mean-field model of
Kuramoto.  This implies that for the system~\eqref{eqn:phase} to
exhibit a phase transition, the interaction function $f$ cannot
be diagonalized.

There is an alternative rigorous approach to Kuramoto's mean-field
model, in which the partial differential equation for the density of
oscillators with certain frequency, which is obtained by taking the
continuum limit $N \to \infty$, is studied to analyze the stability of
the solutions.  See~\cite{strogatz00} for an excellent review in this
direction.

The approach taken here is similar to that
in~\cite[p.~80]{hoppensteadt97b}.  However, some conclusions made
there might be misleading or lack detailed analysis.  This paper is
intended to correct and clarify those points.

The rest of the paper is organized as following.  In
\S\ref{sec:sep}, we introduce an appropriate change of variables
to separate a time-like variable from the rest of the system.  In
\S\ref{sec:diag}, we define diagonalizable interaction and show
how complete separation of variables can be achieved.  We also
establish some properties of diagonalizable systems.  Then, in
\S\ref{sec:rand}, we introduce randomness of the natural
frequencies of the oscillators and state our main results.  Finally,
we discuss some approximate behavior of the system for large $N$ in
\S\ref{sec:dis}, and \S\ref{sec:conc} is reserved for
concluding remarks.

\section{Separation of the Time-like Variable}
\label{sec:sep}

In this and the following sections, we consider the
system~\eqref{eqn:phase} of $N$ phase oscillators, where the natural
frequency vector $\omega$ and
the coupling strength $\eps$ are fixed (nonrandom)
constants.  We will consider $\omega$ to be a random vector in
\S\ref{sec:rand} in order to make probabilistic statements about
the system.

Let us suppose that the interaction function $f$ satisfies two conditions,
\begin{itemize}
\item[(C1)] $\vec{1}^T f(\theta) = \vec{0}$ for all $\theta \in \T^N$ and 
\item[(C2)] $f(v \vec{1} + \theta) = f(\theta)$ for all $\theta \in \T^N$,
\end{itemize}
where $\vec{1} = (1/\sqrt{N}, \ldots, 1/\sqrt{N})^T$.  The condition
(C1) says that the interaction function is orthogonal to the vector
$\vec{1}$.  The second condition (C2) expresses the translation
invariance of $f$ along the direction of $\vec{1}$.  If, for example,
the interaction function has the form~\eqref{eqn:kuramoto}, these
conditions are satisfied if the functions $h_{ij}$ are odd.  In
particular, the mean-field model of Kuramoto does satisfy these
conditions.

Under conditions (C1) and (C2), the system~\eqref{eqn:phase}
can be separated into two independent systems---one for the
time-like variable and the other for the phase deviations.

Let $W$ be an $N \times (N-1)$ matrix whose columns, denoted by $W_j$,
$j = 1,\ldots,N-1$, form an orthonormal basis of the subspace
$\vec{1}^\bot \equiv \{x \in \R^N : \vec{1}^T x = 0 \}$.  In other
words, $W$ is an $N \times (N-1)$ matrix that satisfies $\vec{1}^T W =
\vec{0}$ and $W^T W = I_{N-1}$, where $I_{N-1}$ is the $(N-1) \times
(N-1)$ identity matrix.  Then, the change of variable
\begin{equation}
\label{eqn:change_of_var}
\theta = v \vec{1} + W u
\end{equation}
converts the system \eqref{eqn:phase} into two systems,
\begin{eqnarray}
\dot{v} &=& \vec{1}^T \omega, \label{eqn:v}\\
\dot{u} &=& W^T \omega + \eps W^T f(W u) \label{eqn:u},
\end{eqnarray}
which can be solved separately.  

Systems satisfying the conditions (C1) and (C2) arise in mathematical
neuroscience~\cite{hoppensteadt97b,hoppensteadt97a}, in which $\theta$
often takes the form $\omega t + \phi$ in the limit $t \to \infty$,
where $\omega$ is the vector of carrier frequencies and $\phi$ is the
vector of phase deviations.  The equation~\eqref{eqn:u}, in some
sense, governs the behavior of the phase deviations.

The solution to~\eqref{eqn:v} is $v(t) = v(0) +
(\vec{1}^T \omega) t$, and hence the variable $v$ is time-like if
$\vec{1}^T \omega \neq 0$ or, equivalently, if the average natural
frequency $\sum_i \omega_i/N$ is nonzero.  Thus, the behavior of the
solution of~\eqref{eqn:phase} is essentially determined
by~\eqref{eqn:u}.

Recall that the solution is called coherent if $\lim_{t \to \infty}
\theta_i(t)/\theta_j(t) = 1$.  This can be rephrased in terms of the
vector $\mu \equiv \lim_{t \to \infty} u(t)/t$ of output frequencies
of the $u$-equation~\eqref{eqn:u}, if it exists.  

\begin{lemma} \label{lem:coherence}
Let $u(t)$ be a solution of~\eqref{eqn:u}, and suppose that $\mu \equiv
\lim_{t \to \infty} u(t)/t$ exists.  Then, the solution of
~\eqref{eqn:phase} is coherent if and only if $\mu = 0$.
\end{lemma}
\begin{proof}
Let $\Omega = \lim_{t \to \infty} \theta(t)/t = (\vec{1}^T \omega)
\vec{1} + W \mu$.  Then $\mu = 0$ implies that $\Omega = (\vec{1}^T
\omega) \vec{1}$, which in turn implies that $\lim_{t \to \infty}
\theta_i(t)/\theta_j(t) = 1$.

Conversely, if $\lim_{t \to \infty} \theta_i(t)/\theta_j(t) = 1$,
$\Omega$ must be a multiple of $\vec{1}$.  Since $W$ is orthogonal to
$\vec{1}$, this implies that $W \mu = 0$.  Since $W$ is invertible,
it follows that $\mu = 0$.
\end{proof}

As an immediate consequence of Lemma~\ref{lem:coherence}, if the
solution of~\eqref{eqn:u} tends to an equilibrium, then the
corresponding solution of~\eqref{eqn:phase} is coherent.  For example,
if the interaction function $f$ that satisfies (C1) and (C2) is a
gradient vector field, i.e., $f(\theta) = - \nabla V_0(\theta)$ for
some potential function $V_0: \T^N \rightarrow \R$, then the
$u$-equation~\eqref{eqn:u} is also a gradient system:
\begin{equation}
\label{eqn:gradient_sys}
\dot{u} = - \nabla \left[ - \omega^T W u + \eps V(u) \right],
\end{equation}
where $V(u) = V_0(W u)$.  A minimum $u^*$ of the potential function $-
\omega^T W u + \eps V(u)$ then corresponds to the vector of phase
deviations for a coherent solution of the original oscillator
system~\eqref{eqn:phase}.  As in the proof of
Lemma~\ref{lem:coherence}, we have $\Omega = (\vec{1}^T \omega)
\vec{1}$ in this case, meaning that the output frequency of every
oscillator tends to the mean natural frequency $\bar{\omega} \equiv
\sum_i \omega_i/N$ of the oscillators.  For the interaction of the
form~\eqref{eqn:kuramoto}, the potential function takes the form
\begin{equation*}
V(u) 
= - \frac{1}{2N} \sum_{i=1}^N \sum_{j=1}^N 
  H_{ij} (\theta_i - \theta_j)
= - \frac{1}{2N} \sum_{i=1}^N \sum_{j=1}^N 
  H_{ij} \left(\sum_{k=1}^{N-1} [W_{ik} - W_{jk}] u_k \right),
\end{equation*}
where $H_{ij}(x) \equiv \int_0^x h_{ij}(y) dy$.

\section{Diagonalizable interaction}
\label{sec:diag}

Assuming that the interactions among the oscillators are diagonalizable
enables us to carry out a rigorous analysis of the system.
\begin{definition}
We say that the system~\eqref{eqn:phase} (or the interaction function
$f$) is \emph{diagonalizable} if there exist an $N \times (N-1)$
matrix $W$ and real, continuous, periodic functions $p_j$ such that
\begin{romannum}
\item $\vec{1}^T W = \vec{0}$, 
\item $W^T W = I_{N-1}$, and
\item $f(W u) = W p(u)$ with $p(u) = (p_1(u_1), \ldots, p_{N-1}(u_{N-1}))^T$.
\end{romannum}
\end{definition}
\noindent
For example, $W = W^{(N)}$ defined by
\begin{align}\label{eqn:example}
W^{(N)}_{jk} &= \frac{1}{\sqrt{N}}\left(\sin\frac{2\pi jk}{N} + \cos\frac{2\pi jk}{N}\right)\nonumber\\
             &= \frac{2}{\sqrt{N}}\sin\left(\frac{2\pi jk}{N} + \frac{\pi}{4}\right)
\end{align} 
satisfies these conditions.

When the system~\eqref{eqn:phase} is diagonalizable, the equations for
the components of $u$ become independent of other components:
\begin{equation}
\label{eqn:uj}
\dot{u_j} = a_j + \eps p_j(u_j), \quad j = 1, \ldots, N-1,
\end{equation}
where we set $a_j = W_j^T \omega$.  Thus, the problem is reduced to
solving a scalar differential equation for each $j$.  The following
lemma applies to each equation in~\eqref{eqn:uj}.

\begin{lemma}
\label{lem:ode}
Let $\eps > 0$, and let $a$ be a real number.  Let $p(u)$ be a real,
continuous, periodic function with period $L > 0$.  Define
$\displaystyle m = \min_{0 \le u < L} p(u)$, $\displaystyle M =
\max_{0 \le u < L} p(u)$.  For any solution $u(t)$ of $\dot{u} = a +
\eps p(u)$, the limit $\mu_p(a, \eps) \equiv \lim_{t \to \infty}
u(t)/t$ exists and
\begin{equation}
\mu_p(a, \eps) = \begin{cases}
	L/T(a, \eps), & a < - \eps M, \; a > - \eps m,\\
	0,        & - \eps M \le a \le - \eps m,
\end{cases}
\end{equation}
where
\[
T(a, \eps) \equiv \int_0^L \frac{du}{a + \eps p(u)}
\]
is the ``period'' of the solution in the case of $a < - \eps M$ or $a
> - \eps m$, in the sense that $u(t + T(a, \eps)) = u(t) + L$.
\end{lemma}

\begin{proof}
If $- \eps M \le a \le - \eps m$, then any solution $u(t)$ tends to a
zero of the function $a + \eps p(u)$.  Hence, $\mu_p(a, \eps) = 0$.

For notational simplicity, let us drop the dependence of $T(a, \eps)$
on $a$ and $\eps$ below.  Suppose $a > - \eps m$, so that $a + \eps
p(u) > 0$ for all $u$.  It is straightforward to show that the
function $u(t)$ defined implicitly by the formula
\[
 \int_{u_0}^{u(t)} \frac{du}{a + \eps p(u)} = t
\]
is the unique solution of $\dot{u} = a + \eps p(u)$ with the initial
condition $u(0) = u_0$, and that it satisfies $u(t + T) = u(t) + L$.
We have
\begin{eqnarray*}
\mu_p(a, \eps) &=& \lim_{t \to \infty} u(t)/t\\
  &=& \lim_{t \to \infty} \frac{1}{t}
        \left(u_0 + \int_0^t [a + \eps p(u(s))] ds \right) \\
  &=& a + \eps \lim_{t \to \infty} \frac{1}{t} \int_0^t p(u(s)) ds.
\end{eqnarray*}
Let $n$ be the largest integer for which $nT \le t$.  Then, by
changing the variables in each integral using the translation by
multiples of $T$, we see that
\begin{eqnarray*}
\int_0^t p(u(s)) ds 
  &=& \sum_{k = 1}^n \int_{(k-1)T}^{kT} p(u(s)) ds
      + \int_{nT}^{t} p(u(s)) ds \\
  &=& n \int_{0}^{T} p(u(s)) ds + \int_{0}^{t-nT} p(u(s)) ds.
\end{eqnarray*}
Consequently,
\begin{eqnarray*}
\lefteqn{\left| \frac{1}{t} \int_0^t p(u(s)) ds 
  - \frac{1}{T} \int_0^T p(u(s)) ds \right|} \\
    &=&   \left| \left( \frac{n}{t} - \frac{1}{T} \right) \int_0^T p(u(s)) ds 
	       + \frac{1}{t} \int_0^{t - nT} p(u(s)) ds \right| \\
    &\le& \frac{|nT - t|}{tT} \int_0^T |p(u(s))| ds 
               + \frac{1}{t} \int_0^{t - nT} |p(u(s))| ds\\
    &\le& \frac{T}{tT}\ T \max\{|m|,|M|\} + \frac{1}{t}\ T  \max\{|m|,|M|\} \\
    &\to& 0,
\end{eqnarray*}
as $t \to \infty$, showing that the limit $\lim_{t \to \infty}
\frac{1}{t} \int_0^t p(u(s)) ds$ exists and is equal to $\frac{1}{T}
\int_0^T p(u(s)) ds$.  Thus, by changing variables from $s$ to $u$ and
translating by $u_0$, we see that $\mu_p(a, \eps)$ exists and
\begin{eqnarray*}
\mu_p(a, \eps) 
	&=& a + \frac{\eps}{T} \int_0^T p(u(s)) ds \\
	&=& a + \frac{\eps}{T} \int_0^L \frac{p(u)du}{a + \eps p(u)} \\
	&=& L/T.
\end{eqnarray*}

If $a < - \eps M$, then, by replacing $u$ with $-u$, $a$ with $-a$,
$\eps$ with $- \eps$, and $m$ with $M$, the problem reduces to the
previous case.  The lemma is proved.
\end{proof}

\newcommand{\sgn}{\mathrm{sgn}}  

\begin{figure}
\setlength{\unitlength}{0.49\textwidth}

\begin{picture}(1,1)
\put(0.5,0.5){\makebox(0,0){
\resizebox{0.9\unitlength}{!}{
\includegraphics{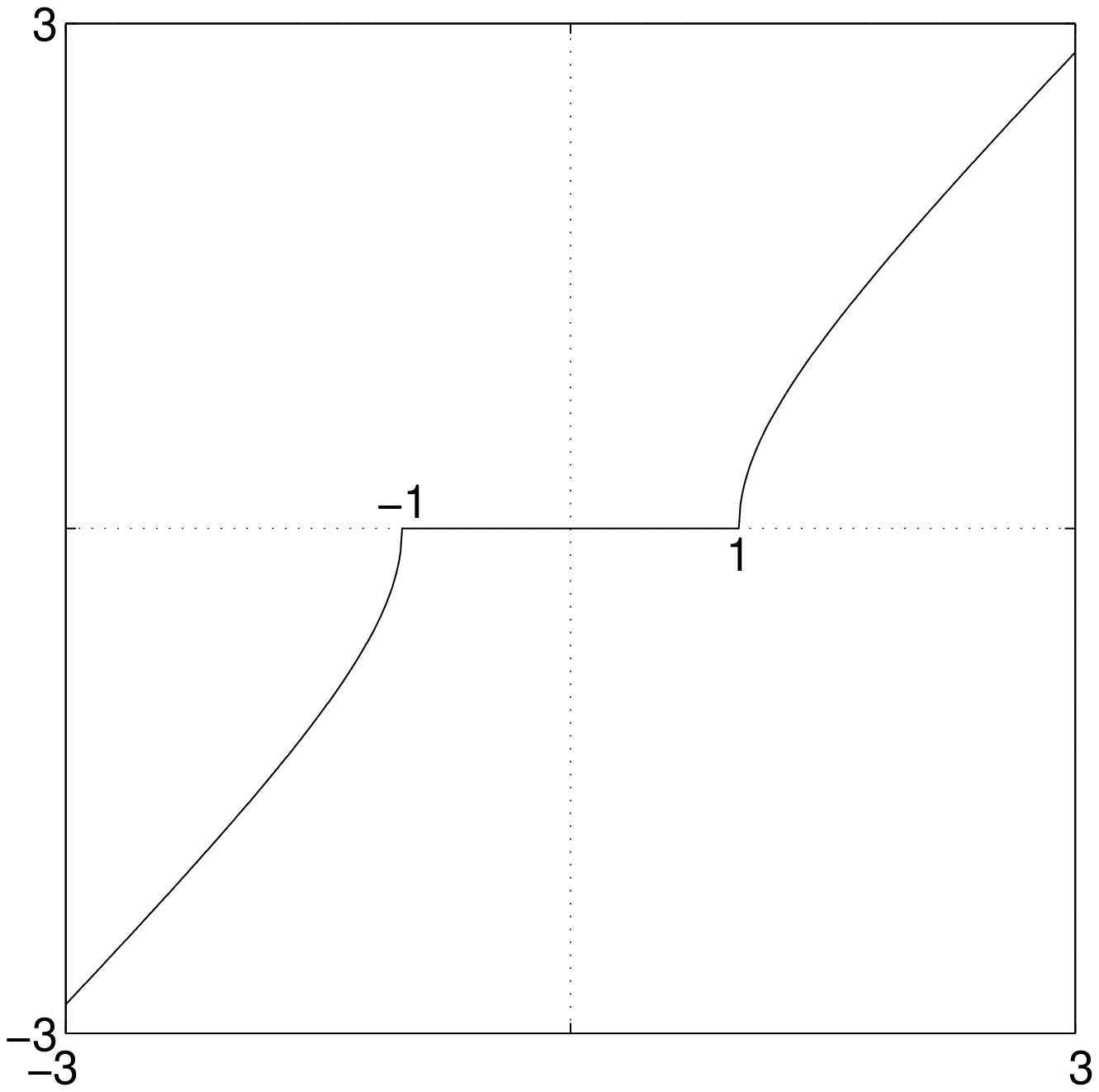}}}}
\put(0.525,1){\makebox(0,0)[t]{(a)}}
\put(0.53,0.02){\makebox(0,0)[b]{$a$}}
\put(0,0.5){\makebox(0,0)[l]{\rotatebox{90}{$\mu_p(a, 1)$}}}
\end{picture}
\begin{picture}(1,1)
\put(0.5,0.5){\makebox(0,0){
\resizebox{0.9\unitlength}{!}{
\includegraphics{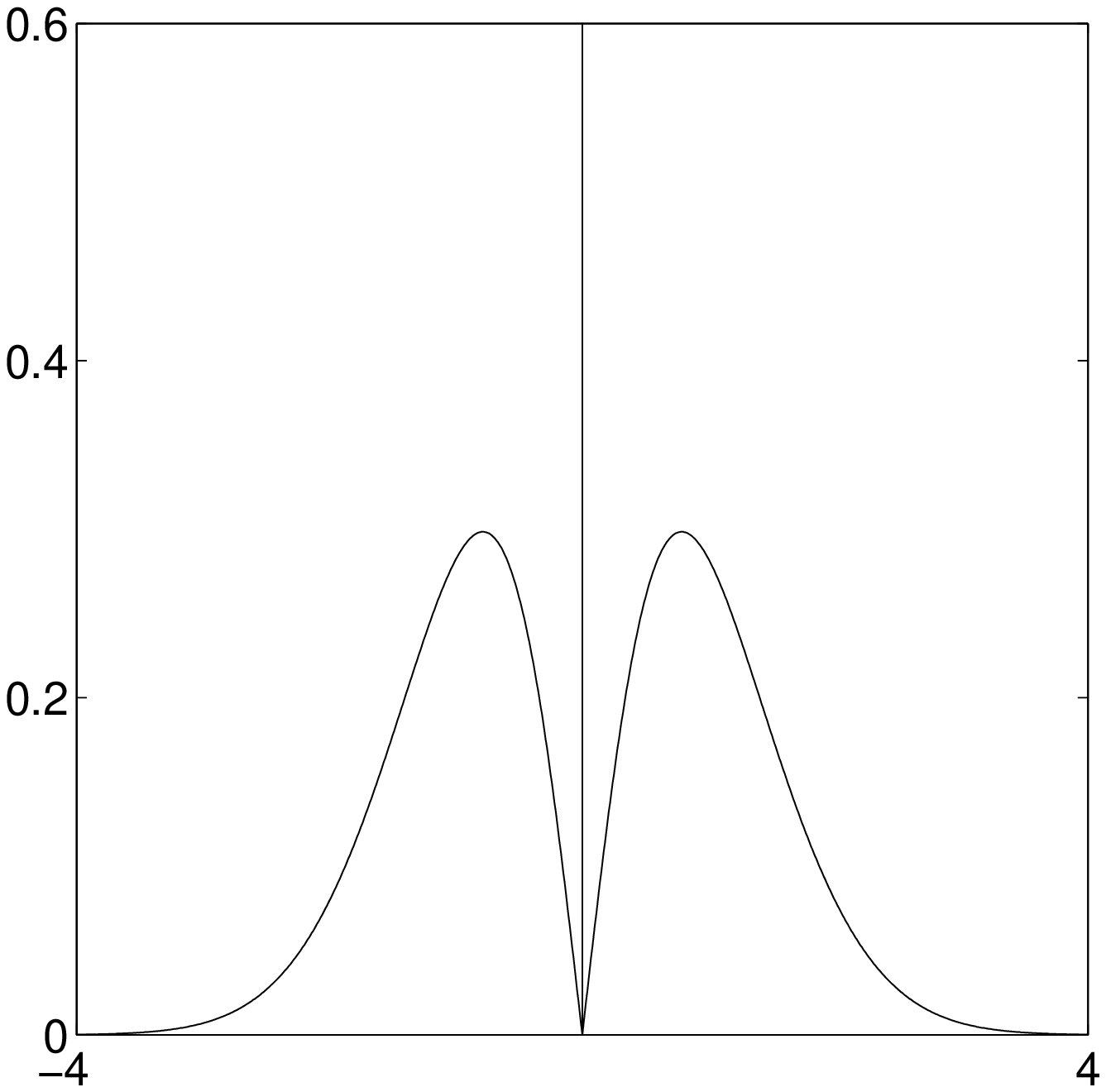}}}}
\put(0.53,1){\makebox(0,0)[t]{(b)}}
\put(0.53,0.02){\makebox(0,0)[b]{$\mu$}}
\put(0,0.5){\makebox(0,0)[l]{\rotatebox{90}{$g(\mu; 1)$}}}
\end{picture}

\caption{(a) The graph of the input-output frequency function
$\mu_p(a, \eps)$ vs $a$ for $\eps = 1$ and $p(u) = \sin(u)$.  (b) The
corresponding density $g(\mu; 1)$ when $a$ is the standard Gaussian
random variable.}
\label{fig:mu}
\end{figure}

The function $\mu_p(\cdot, \eps)$ in Lemma~\ref{lem:ode}, which can
easily be shown to be differentiable with positive derivative outside
the interval $[-\eps M, \eps M]$, determines the relationship between
the input frequency $a$ and the output frequency $\mu_p(a,\eps)$.  If
we take $p(u) = \sin(u)$, for example, the integration in the
expression of $\mu_p$ can be carried out, and we get
\[
\mu_p(a, \eps) = \begin{cases}
	- \sqrt{a^2 - \eps^2},   & a < - \eps, \\
	0,                          & -\eps \le a < \eps, \\
	\sqrt{a^2 - \eps^2}, & a \ge \eps,
\end{cases}
\]
the graph of which is given in Fig.~\ref{fig:mu}(a) for $\eps = 1$.  With
this function, the output frequency vector $\mu$ of the
$u$-equation~\eqref{eqn:u} can be written as
\[ \mu = \mu(a) =
(\mu_{p_1}(a_1, \eps), \ldots, \mu_{p_{N-1}}(a_{N-1}, \eps))^T. \]

It is important to note here that $\mu$ does \emph{not} depend on the
initial condition $u(0) = u_0$, which implies that it is also
independent of the initial condition for $\theta$.  In other words,
the initial condition for the system~\eqref{eqn:phase} does not affect
the behavior of its solution, as far as its coherence properties are
concerned.  Therefore, in this sense, coherence, partial coherence and
incoherence are properties of the system rather than of individual
solutions for a diagonalizable system.

\section{Randomly Distributed Frequencies}
\label{sec:rand}

In this section we consider $\omega$ to be a random vector in $\R^N$.
We take the components $\omega_1, \ldots, \omega_N$ of $\omega$ to be
independent and identically distributed (i.i.d.) random variables with mean 0 and variance $\sigma^2 > 0$.  In the general case
of mean $\omega_0 \neq 0$, the problem can always be reduced to the
zero-mean case by the translation of $\theta$ by $- \omega_0 t$.

Since $\omega$ is random, the vectors $a$ and $\mu$ are also random
vectors in $\R^{N-1}$.  Lemma \ref{lem:ode} along with the relation $a
= W^T \omega$ can be used to determine the distribution of $\mu$ from
the distribution of $\omega$.  For example, if each $\omega_j$ is
standard Gaussian, then so is each $a_j$, in which case the density $g(\mu;
\eps)$ for the random variable $\mu_p(a_j, 1)$ when $p(u) = \sin(u)$
can be computed.  The result is
\[
g(\mu; \eps) = 
  \frac{|\mu|e^{- (\mu^2 + \eps^2)/2}}{\sqrt{2 \pi(\mu^2 + \eps^2)}}
  + \delta(\mu) \, \erf \! \left( \frac{\eps}{\sqrt{2}} \right),
\]
where $\delta(\mu)$ is Dirac's delta function.  The graph of this
density is shown in Fig.~\ref{fig:mu} for $\eps = 1$.

Our main goal is this section is to compute the probabilities that the
system~\eqref{eqn:phase} is coherent, partially coherent, or
incoherent.  The following theorem reveals a curious property of a
generic diagonalizable system of phase oscillators.

\begin{theorem}
\label{thm:1}
Let the natural frequency vector $\omega$ be a random vector in
$\R^N$, whose components are i.i.d.\ with a common continuous
distribution.  Suppose that~\eqref{eqn:phase} is a diagonalizable system of
$N$ phase oscillators such that $W$ satisfies the condition that
$W_{ki} \neq W_{kj}$ for all $k = 1,2,\ldots,N$ and for all $i,j =
1,2,\ldots,N-1$ such that $i \neq j$.  Then, the partial coherence of
the system almost surely implies coherence; i.e., given that the
system is partially coherent, the probability that it is coherent is
one.
\end{theorem}

\begin{proof}
Once again, let $\Omega(\omega) = \lim_{t \to \infty} \theta(t)/t =
(\vec{1}^T \omega) \vec{1} + W \mu(\omega)$.  Let $S_c$ be the set of
$\omega$ in $\R^N$ that corresponds to coherent systems, i.e., $S_c =
\bigl\{ \omega \in R^N : \Omega(\omega) = \vec{0} \bigr\}$.  By
Lemma~\ref{lem:coherence}, we may also write $S_c = \bigl\{ \omega \in R^N
: \mu(\omega) = \vec{0} \bigr\}$.  Let $S_{pc}$ be the set corresponding to
partially coherent systems, that is, $S_{pc} = \bigl\{ \omega \in R^N
: \Omega_i(\omega) = \Omega_j(\omega)$ for some $i \neq j
\bigr\}$.  It is easy to see that we can also rewrite this in terms
of $\mu$ as
\begin{eqnarray*}
S_{pc} 
  &=& \left\{ \omega \in R^N : \text{ There are }i \neq j \text{ s.t. }
         \sum_{k=1}^{N-1} (W_{ki} - W_{kj}) \mu_k(\omega) = 0 \right\} \\
  &=& \bigcup_{i \neq j} \left\{ \omega \in R^N : 
         \sum_{k=1}^{N-1} (W_{ki} - W_{kj}) \mu_k(\omega) = 0 \right\}
  \equiv \bigcup_{i \neq j} S_{pc}^{(i,j)}.
\end{eqnarray*}

The probability that the system is coherent, given that the system is
partially coherent, is $P(S_c)/P(S_{pc})$ since $S_c \subset S_{pc}$.
This probability is one if and only if $P(S_{pc} \setminus S_c) = 0$,
which would be satisfied if $P(S_{pc}^{(i,j)} \setminus S_c) = 0$ for
every pair $i \neq j$.  We shall show this next.

Let us fix $i$ and $j$.  For any $A \subset \{ 1,2,\ldots,N-1 \}$,
denote by $Z_k$ the subspace $\bigl\{ \mu \in \R^{N-1} : \mu_k = 0 \bigr\}$, and
let $Z_A = \bigcup_{k \in A} Z_k$ and $Z'_A = \bigcap_{k \notin A}
Z_k$.  Let $R_A$ denote the subspace $\bigl\{ \mu \in \R^{N-1} : \sum_{k
\in A} (W_{ki} - W_{kj}) \mu_k = 0 \bigr\}$.  Define $Q_A = R_A \cap Z'_A
\setminus \zero$.  We will show that $P(Q_A) = 0$ for any choice of
$A$.  $P(S_{pc}^{(i,j)} \setminus S_c) = 0$ follows from this by
taking $A = \{ 1,2,\ldots,N-1 \}$.

We shall prove $P(Q_A) = 0$ by induction on $n = |A|$, the cardinality
of $A$.  Suppose first that $n = 1$ and, say, $A = \{ 1 \}$.  Since
$W_{1 i} - W_{1 j} \neq 0$, we have $R_A = \bigl\{ \mu \in \R^{N-1} : \mu_1
= 0 \bigr\} = Z_1$ and $Z'_A = \bigcap_{k=2}^{N-1} Z_k$.  Thus, $Q_A =
\bigcap_{k=1}^{N-1} Z_k \setminus \zero = \varnothing$, which implies
$P(Q_A) = 0$.  The same holds for any other $A$ with $|A| = 1$.

Suppose that $P(Q_A) = 0$ for any $A$ with $|A| = n-1$, and consider the
case $|A| = n$.  We have
\begin{eqnarray*}
P(Q_A) &=& P(Q_A \cap Z_A) + P(Q_A \setminus Z_A) \\
       &=& P\left( \bigcup_{k \in A} Q_A \cap Z_k \right) 
            + P(Q_A \setminus Z_A) \\
       &\le& \sum_{k \in A} P(Q_A \cap Z_k) + P(Q_A \setminus Z_A).
\end{eqnarray*}
We see that $P(Q_A \cap Z_k) = 0$ for each $k \in A$, by the induction
hypothesis, since we can write $Q_A \cap Z_k = R_A \cap Z'_A \cap Z_k
\setminus \zero = R_{A_k} \cap Z'_{A_k} \setminus \zero = Q_{A_k}$
with $A_k = A \setminus \{ k \}$, for which we have $|A_k| = n-1$.
Thus, if we can show $P(Q_A \setminus Z_A) = 0$, then we are done.

We show $P(Q_A \setminus Z_A) = 0$ in three steps.  First, since
$Z'_A$ is an $n$-dimensional subspace and $Q_A \subset R_A \cap Z'_A$
is an $(n-1)$-dimensional subspace, the $n$-dimensional Lebesgue
measure of $Q_A$ in $Z'_A$ must be zero.

Next, note that the conditional probability distribution of $\mu$,
given $\mu \in Z'_A$, is continuous with respect to the Lebesgue
measure outside the set $Z_A$.  This can be seen by noting the
following: (1) each component can be written as $\mu_j =
\mu_{p_j}(a_j)$ by Lemma~\ref{lem:ode}, (2) $\mu_{p_j}^{-1}$ exists
and is differentiable except at the origin, again by
Lemma~\ref{lem:ode}, and (3) the conditional distribution of $a_j =
\sum_k W_{kj} \omega_k$, given that $\mu(\omega) \in Z'_A$ (which is
equivalent to $-\eps M_k \le \sum_l W_{lk} \omega_l \le -\eps m_k$ for
all $k \notin A$), is continuous everywhere.

Finally, combining these two observations, we see that $P(Q_A
\setminus Z_A \;|\; Z'_A) = 0$, which implies that $P(Q_A \setminus
Z_A) = P(Z'_A) P(Q_A \setminus Z_A \;|\; Z'_A) = 0$.  This completes the proof of the theorem.
\end{proof}

\newcommand{\cS}{{\mathcal S}}

We next describe the behavior of a generic diagonalizable
system in the limit of $N \to \infty$.  In order to formalize the
process of taking the limit, we need to choose a sequence of systems
of the form~\eqref{eqn:phase}.  Such a sequence can be characterized
by the following.

\begin{enumerate}
\item Consider a sequence $\{ W^{(N)} \}_{N = 1,2,\ldots}$ of matrices
with the following properties:

  \begin{enumerate}

    \item Each $W^{(N)}$ is an $N \times (N-1)$ matrix with orthonormal columns.

    \item $\vec{1}_N^T W^{(N)} = \vec{0}$ for all $N$.

    \item Each $W^{(N)}$ satisfies the condition for $W$ in
    	  Theorem~\ref{thm:1}.

    \item $||W^{(N)}||_\infty \to 0$ as $N \to \infty$.
	  (Here $||\cdot||_\infty$ denotes the maximum matrix norm defined
	  by $||A||_\infty = \max |A_{ij}|$, where the maximum is
	  taken over all elements of $A$.)  This is like a mixing condition that will be necessary
	  later in order to apply Proposition~\ref{prop:clt}.

  \end{enumerate}

\item Consider a sequence $\{ p_j \}_{j=1,2,\ldots}$ of real, continuous,
periodic functions such that the corresponding sequence of norms
$\displaystyle ||p_j|| \equiv \max_j |p_j(u)|$ is bounded.

\item Consider a sequence $\{ \omega_j \}_{j=1,2,\ldots}$ of i.i.d.\
random variables with mean $\omega_0$ and variance $\sigma^2$.
\end{enumerate}
The sequence of matrices $W^{(N)}$ defined by \eqref{eqn:example} satisfies the conditions above.
Given such sequences, for each $\eps > 0$ and $N$, we define $\cS_{N, \eps}$ to be the diagonalizable system~\eqref{eqn:phase} of phase
oscillators using the natural frequency vector $\omega = (\omega_1,
\ldots, \omega_N)^T$, the functions $\{ p_1, \ldots, p_{N-1} \}$ and
the matrix $W^{(N)}$.  We are now ready to state and prove
our main theorem.

\begin{theorem}
\label{thm:2}
Let $\cS_{N,\eps}$ be defined as above.  Then, for any fixed $\eps > 0$,
$\cS_{N, \eps}$ is almost surely incoherent as $N \to \infty$; i.e.,
the probability that $\cS_{N,\eps}$ is incoherent tends to one in the
limit of $N \to \infty$.
\end{theorem}

\begin{proof}
As mentioned before, we may assume $\omega_0 = 0$ without loss of
generality, since the $\omega_0 \neq 0$ case can always be reduced to
the $\omega_0 = 0$ case.

From Theorem~\ref{thm:1}, we know that the probability that
$\cS_{N,\eps}$ is \emph{not} incoherent is equal to the probability
$q_c$ that it is coherent.  We need to show that $q_c \to 0$ as $N \to
\infty$.

Let $N_0 < N$ be fixed.  From Lemma~\ref{lem:coherence}, it follows
that $q_c = P(\mu = \vec{0})$.  Since the sequence $\{ ||p_j||
\}_{j=1,2,\ldots}$ is bounded, we can define $M = \sup M_j$ and $m =
\inf m_j$, where $\displaystyle m_j = \min_{0 \le u < L} p_j(u)$,
$\displaystyle M_j = \max_{0 \le u < L} p_j(u)$ for each $j$.  Then,
Lemma~\ref{lem:ode} implies
\begin{eqnarray*}
q_c &=& P(\mu = \vec{0})\\
    &=& P( - \eps M_j \le a_j^{(N)} \le - \eps m_j,\ j = 1, \ldots, N-1)\\
    &\le& P( - \eps M \le a_j^{(N)} \le - \eps m,\ j = 1, \ldots, N_0),
\end{eqnarray*} 
where $a_j^{(N)} = \left( W^{(N)}_j \right)^T \omega$.

The following Proposition shows that for each $j$, $a_j^{(N)}$
converges to a Gaussian random variable with mean 0 and variance
$\sigma^2$.

\begin{proposition}
\label{prop:clt}
Let $X_1, X_2, \ldots$ be a sequence of i.i.d.\ random variables with
$EX_j = 0$ and $\var{X_j} = E(X_j^2) = \sigma^2$.  Suppose that, for
each $N$, real numbers $b_{N,1}, \ldots, b_{N,N}$ satisfy
$\sum_{j=1}^N b_{N,j}^2 = 1$.  Also, suppose that 
\[ \lim_{N \to \infty} \max_{1 \le j \le N} |b_{N,j}| = 0. \]  
Then we have
\[ S_N = \sum_{j=1}^N
b_{N,j} X_j \xrightarrow{d} {\mathcal N}(0, \sigma^2) \] 
as $N \to \infty$.
\end{proposition}
\begin{proof}
Let $Y_{N,j} = b_{N,j} X_j$.  We will apply the Lindeberg--Feller
central limit theorem (see, for example, \cite[p.~98]{durrett}) to $Y_{N,j}$.
For this we need to check three conditions.  The first is
$EY_{N,j} = b_{N,j} EX_j = 0$.  The second condition is satisfied
because $\sum_{j=1}^N EY_{N,j}^2 = \sum_{j=1}^N b_{N,j}^2 EX_j^2 =
\sigma^2 > 0$.  To show that the third condition is satisfied, let
$\eps > 0$ be fixed.
We have
\[
\sum_{j=1}^N E\Bigl(|Y_{N,j}|^2 \Big|\; |Y_{N,j}| > \eps \Bigr)
  = \sum_{j=1}^N b_{N,j}^2 E\left( |X_j|^2 \bigg|\; 
                |X_j| > \frac{\eps}{|b_{N,j}|} \right),
\]
where $E(X|A)$ denotes the conditional expectation of $X$, given $A$.
Let $j$ be fixed.  For each $N$, set $Z_N = |X_j|^2$ if $|X_j| >
\eps/|b_{N,j}|$, and 0 otherwise.  Since $|b_{N,j}| \to 0$, $Z_N \le
|X_j|^2$ for each $N$, and $Z_N \to 0$ almost surely, we may use the
dominated convergence theorem to show that for each $j = 1,2,\ldots$,
$E Z_N = E\bigl(|X_j|^2 \;\big|\; |X_j| > \eps/|b_{N,j}|\bigr) \to 0$ as $N \to
\infty$.  Thus, the third condition $\sum_{j=1}^N E\bigl(|Y_{N,j}|^2 \;\big|\;
|Y_{N,j}| > \eps\bigr) \to 0$ is satisfied.  The conclusion now follows
directly from application of the Lindeberg-Feller theorem.
\end{proof}

For each $j = 1, \dots, N-1$, we take $b_{N,i} =
W^{(N)}_{ij}$ and $X_i = \omega_i$ in Proposition~\ref{prop:clt}, and we see that $a_j^{(N)}$ converges in distribution to $a_j^{(\infty)}$ as
$N \to \infty$, where $a_j^{(\infty)}$ is a Gaussian random variable.
Moreover, due to the orthogonality of $W^{(N)}$, $a_1^{(N)}, \ldots,
a_{N-1}^{(N)}$, in some sense, become independent in the limit.

\begin{lemma}
The random variables $a_1^{(\infty)},a_2^{(\infty)}, \ldots$ are
independent.
\end{lemma}
\begin{proof}
We need to show that for any finite $A \subset \N$, the collection $\{
a_k^{(\infty)} \}_{k \in A}$ is a set of independent random variables.
For simplicity, we only prove this for $A = \{1,2\}$, but a similar argument
works for a general case. 

Let $t_1$ and $t_2$ be given.  Set
\[
b_{N,j} = \frac{t_1 W_{1j}^{(N)} + t_2 W_{2j}^{(N)}}{\sqrt{t_1^2 + t_2^2}}.
\]
Then, as $N \to \infty$, $\max_j |b_{N,j}|$ approaches zero because
$\max_j |W_{1j}|$ and $\max_j |W_{2j}|$ go to zero.  Also, it is easy
to check that $\sum_{j=1}^N b_{N,j}^2 = 1$ for all $N$.  Applying
Proposition~\ref{prop:clt}, we see that
\[
\frac{t_1 a_1^{(N)} + t_2 a_2^{(N)}}{\sqrt{t_1^2 + t_2^2}}
 = \sum_{j=1}^N b_{N,j} \omega_j
 \xrightarrow{d} {\mathcal N}(0, \sigma^2),
\]
which implies that $t_1 a_1^{(N)} + t_2 a_2^{(N)} \xrightarrow{d}
{\mathcal N}(0, \sigma^2(t_1^2 + t_2^2))$, which in turn implies the
convergence of the joint characteristic function of $a_1^{(N)}$ and
$a_2^{(N)}$ as $N \to \infty$.  Specifically,
\[
E e^{i t_1 a_1^{(\infty)} + i t_2 a_2^{(\infty)}}
  = \lim_{N \to \infty} E e^{i(t_1 a_1^{(N)} + t_2 a_2^{(N)})} 
  = e^{- \sigma^2(t_1^2 + t_2^2)/2}
  = e^{- \sigma^2 t_1^2/2} \, e^{- \sigma^2 t_2^2/2}.
\]
Therefore, $a_1^{(\infty)}$ and $a_2^{(\infty)}$ are independent.
\end{proof}

Let us come back to the proof of Theorem~\ref{thm:2}.  As a
consequence of $a_1^{(\infty)},a_2^{(\infty)},\ldots$ being
independent Gaussian random variables, we have
\begin{eqnarray*}
q_c &\le& P(\{ - \eps M \le a_j^{(N)} \le - \eps m,\ j = 1, \ldots, N_0 \})\\
    &\to& \left[
            \frac{1}{\sqrt{2 \pi \sigma^2}} 
            \int_{- \eps M}^{- \eps m} e^{-x^2/2\sigma^2} dx 
          \right]^{N_0}
\end{eqnarray*} 
as $N \to \infty$.  Since this holds for any fixed $N_0$, and since the
right-hand side goes to zero as $N_0 \to \infty$, we conclude that
$q_c \to 0$.  This completes the proof of Theorem~\ref{thm:2}.
\end{proof}

\section{Example}
\label{sec:example}
A network of voltage-controlled oscillator (VCO) devices can be built as an example of systems with diagonalizable interaction.  The behavior of the $j$th VCO in the network is described by its phase variable $\theta_j$, which satisfies~\cite{hoppensteadt97b,horowitz}:
$$ \dot{\theta_j} = \omega_j + I_j(t), $$
where $\omega_j$ is the center frequency and $I_j(t)$ is the input signal from other VCOs.  The system is diagonalizable if, for example, $I_j(t)$ has the form
$$ I_j(t) = \eps \sum_{k=1}^{N-1} W^{(N)}_{jk} \sin\left(\sum_{\ell=1}^N W^{(N)}_{\ell k} \theta_\ell \right),$$
with $W^{(N)}$ defined by \eqref{eqn:example}.  

This type of interaction can be implemented using commercially available circuit elements, as follows.  The sine terms on the right-hand side can be constructed as sums and products of output voltages:
$$ \sin\left(\sum_{\ell=1}^N W^{(N)}_{\ell k} \theta_\ell \right) = \sum_{b} \prod_{\ell=1}^N \cos\left(W^{(N)}_{\ell k} \theta_\ell - \frac{b_\ell\pi}{2}\right), $$
where the sum is taken over all (ordered) binary $N$-tuples $b=(b_1,\ldots,b_N)$, $b_\ell=0,1$, such that $\sum_\ell b_\ell$ is odd.  This means that $I_j(t)$ is a sum of terms that are products of $\sin(W^{(N)}_{\ell k} \theta_\ell)$ and $\cos(W^{(N)}_{\ell k} \theta_\ell)$.  Signals of the form $\sin(W^{(N)}_{\ell k} \theta_\ell)$ can be obtained by using the amplified version of the input $W^{(N)}_{\ell k}I_\ell(t)$ as the controlling voltage in a separate VCO with center frequency $W^{(N)}_{\ell k}\omega_\ell$.  From these we can get $\cos(W^{(N)}_{\ell k} \theta_\ell)$ by the phase shift of $\pi/2$.  Finally, $I_j(t)$ is obtained by putting these signals through multipliers and adding the outputs.  

\section{Discussion}
\label{sec:dis}

In order to gain additional insights, let us consider the case when the
functions $p_j=p$ do not depend on $j$, and have the range
$[-1,1]$. The arguments in the proof of Theorem~\ref{thm:2} suggest
that for a diagonalizable system~\eqref{eqn:phase} with large $N$, the
probability $q_c$ that it is coherent is approximately
\[
q_c \approx \left[
               \frac{1}{\sqrt{2 \pi \sigma^2}} 
               \int_{-\eps}^{\eps} e^{-x^2/2\sigma^2} dx 
            \right]^{N-1}
    = \left[
          \erf \!\left( \frac{\eps}{\sigma \sqrt{2}} \right)
      \right]^{N-1}
    \equiv \tilde{q}_c(\eps; N),
\]
where $\erf (x)$ is the error function.  Typical graphs of
$\tilde{q}_c(\eps; N)$ are plotted in Fig.~\ref{fig:qc}.  One can see
that for any finite value of $N$, there seems to be a sharp transition
point through which $\tilde{q}_c(\eps; N)$ changes from 0 to 1.
However, unlike the mean-field model of Kuramoto, this point keeps
shifting to the right as $N$ increases, and tends to $\infty$ in the
limit of $N \to \infty$, although it can be shown that this increase
is at most $O(\sqrt{\ln N})$:
\begin{figure}
\setlength{\unitlength}{0.7\textwidth}
\begin{center}
\begin{picture}(1,0.8)
\put(0.5,0.4){\makebox(0,0){
\resizebox{0.9\unitlength}{!}{
\includegraphics{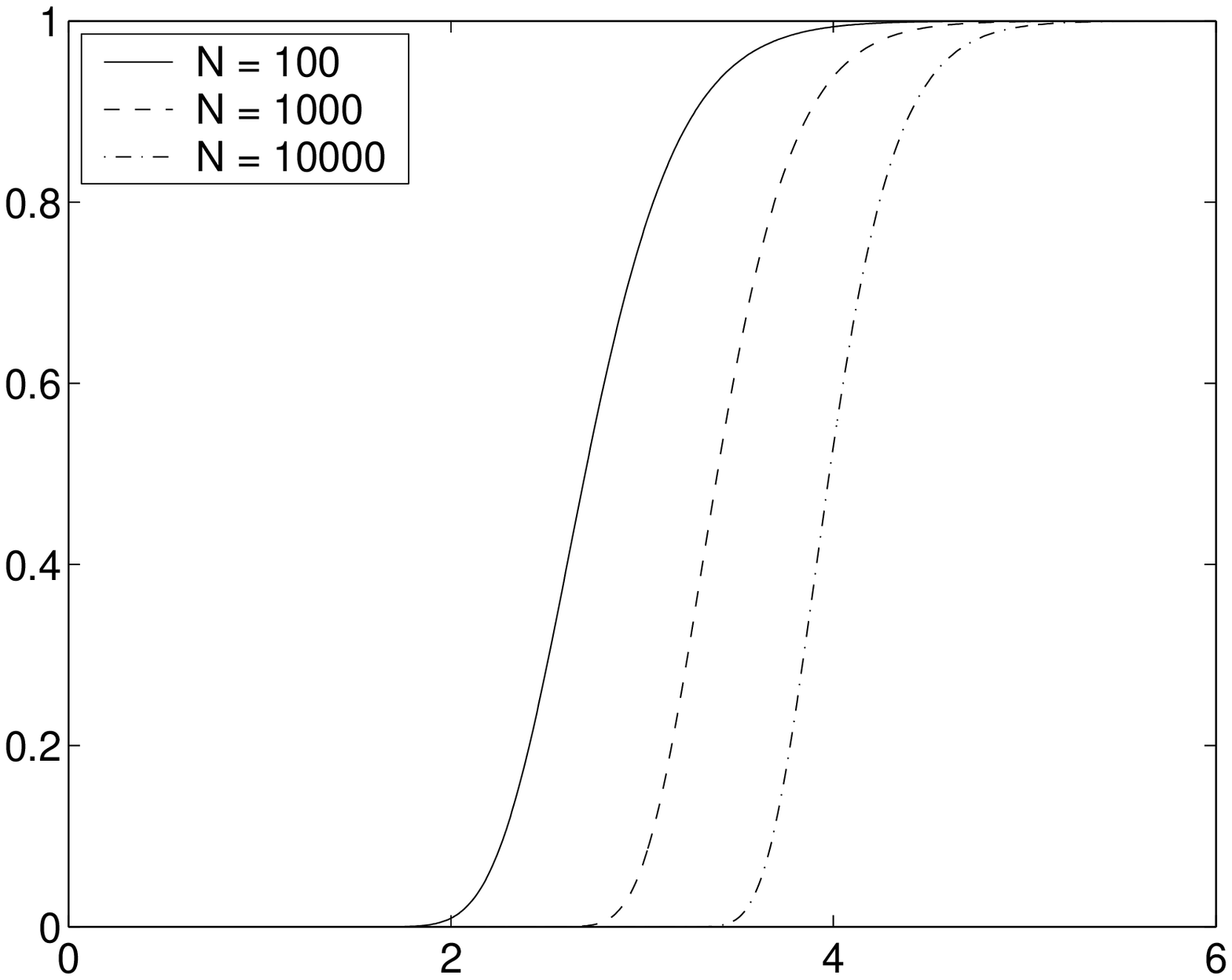}}}}
\put(0.52,0.02){\makebox(0,0)[b]{$\eps$}}
\put(0,0.41){\makebox(0,0)[l]{\rotatebox{90}{$\tilde{q}_c(\eps; N)$}}}
\end{picture}
\end{center}
\caption{Probability of coherence $\tilde{q}_c(\eps; N)$ as a function of $\eps$ for $\sigma = 1$ and $N =
100, 1000, 10000$.}
\label{fig:qc}
\end{figure}
\begin{lemma}
Let $\sigma > 0$ and $0<q<1$ be fixed. Define $\eps_{q,\sigma}(N)$ implicitly by $q =
\tilde{q}_c(\eps_{q,\sigma}(N); N)$.  Then $\eps_{q,\sigma}(N) =
O(\sqrt{\ln N})$ as $N \to \infty$; i.e.,
$\eps_{q,\sigma}(N)/\sqrt{\ln N}$ is bounded as $N \to \infty$.
\end{lemma}
\begin{proof}
Let $x = (\sigma \sqrt{2})^{-1}\eps_{q,\sigma}(N)$.  Then, using a
known estimate for the error function, we have for $x \ge 1$
\[
\erf(x) 
\ge 1 - \frac{2 e^{-x^2}}{\sqrt{\pi}x + \sqrt{\pi x^2 + 4}}
\ge 1 - \frac{2 e^{-x^2}}{\sqrt{\pi} + \sqrt{\pi + 4}}.
\]
Hence, we have the estimate
\[
q \ge (1 - C_0 e^{- x^2})^{N-1} \ge 1 - (N-1)C_0 e^{- x^2},
\]
where
\[
C_0 = \frac{2}{\sqrt{\pi} + \sqrt{\pi + 4}}.
\]
Here we used the relation $(1-x)^n \ge 1 - n x$, which is valid for $n
\ge 0$ and $0 \le x \le 1$.  The estimate for $\eps_{q,\sigma}(N)$ can
be obtained by rearranging:
\[
\eps_{q,\sigma}(N) \le \sigma \sqrt{C_1 + 2 \ln(N-1)},
\]
where $C_1 = 2 \ln C_0 - \ln(1-q)$.  This implies $\eps_{q,\sigma}(N)
= O(\sqrt{\ln N})$.
\end{proof}

\section{Conclusions}
\label{sec:conc}

In this paper, we have defined a class of systems of phase oscillators
characterized by having diagonalizable interactions.  For a system in
this class, complete separation of variables through appropriate
changes of variable is possible, which enables us to draw rigorous
conclusions about the probabilistic properties of the system.  In
particular, we have shown that partial coherence of the system almost
surely implies coherence and, in the limit of large system size, the
system is almost surely incoherent.  A major implication of our result
is that, unlike the mean-field model of Kuramoto, diagonalizable systems
cannot exhibit a sharp transition from incoherence to coherence.  This
provides some insight into what is necessary to see such a transition
in a system of phase oscillators.

\bibliographystyle{siam}
\bibliography{sync,paper}

\end{document}